\documentclass[conference]{IEEEtran}
\IEEEoverridecommandlockouts
\usepackage{cite}
\usepackage{amsmath,amsfonts}
\usepackage{algorithmic}
\usepackage{graphicx}
\usepackage{xcolor}
\usepackage{enumitem}
\def\BibTeX{{\rm B\kern-.05em{\sc i\kern-.025em b}\kern-.08em
    T\kern-.1667em\lower.7ex\hbox{E}\kern-.125emX}}

\usepackage{tikz}
\usetikzlibrary{bayesnet}
\usetikzlibrary{shapes, arrows, calc, positioning,matrix}
\tikzset{
data/.style={circle, draw, text centered, minimum height=3em ,minimum width = .5em, inner sep = 2pt},
empty/.style={circle, text centered, minimum height=3em ,minimum width = .5em, inner sep = 2pt},
}

\usepackage{epsfig}
\usepackage{graphicx}
\usepackage{subfigure}
\usepackage[ruled,vlined]{algorithm2e}
\usepackage{multicol}
\usepackage{multirow}
\usepackage{enumitem}
\usepackage{stfloats}
\usepackage{color}
\usepackage{url}
\usepackage{booktabs}
\usepackage{mathtools}
\usepackage{tkz-euclide}
\usepackage{url}

\usepackage{bm}
\DeclareMathAlphabet\mathbfcal{OMS}{cmsy}{b}{n}
\newcommand{\ten}[1]{\mathbfcal{#1}} 
\newcommand{\mat}[1]{\mathbf{#1}}

\usepackage{pst-plot}

\newcommand{\parm}{{\xi}}
\newcommand{\vecpar}{\boldsymbol{\parm}}

\newcommand{\multiGPC}{\Psi }

\newcommand{\polyInd}{\alpha}
\newcommand{\basisInd}{\boldsymbol{\polyInd}}

\begin{document}
\title{Tensor Methods for Generating\\ Compact Uncertainty Quantification and Deep Learning Models {\footnotesize 
\thanks{$^*$Equally contributed authors. This work was partly supported by NSF CAREER Award No. 1846476, NSF-CCF Awards No. 1763699 and No. 1817037, and an UCSB start-up grant.}
}}

\IEEEspecialpapernotice{(Invited Special Session Paper)}

\author{\IEEEauthorblockN{Chunfeng~Cui$^*$\IEEEauthorrefmark{2},~Cole~Hawkins$^*$\IEEEauthorrefmark{3},~and~Zheng~Zhang\IEEEauthorrefmark{2}}
\IEEEauthorblockA{\IEEEauthorrefmark{2}Department of Electrical and Computer Engineering, 
University of California, 
Santa Barbara, CA 93106}
\IEEEauthorblockA{\IEEEauthorrefmark{3}Department of Mathematics,  University of California, Santa Barbara, CA 93106\\
}
}
\maketitle

\begin{abstract}
Tensor methods have become a promising tool to solve high-dimensional problems in the big data era. By exploiting possible low-rank tensor factorization, many high-dimensional model-based or data-driven problems can be solved to facilitate decision making or machine learning. In this paper, we summarize the recent applications of tensor computation in obtaining compact models for uncertainty quantification and deep learning. In uncertainty analysis where obtaining data samples is expensive, we show how tensor methods can significantly reduce the simulation or measurement cost. To enable the deployment of deep learning on resource-constrained hardware platforms, tensor methods can be used to significantly compress an over-parameterized neural network model or directly train a small-size model from scratch via optimization or statistical techniques. Recent Bayesian tensorized neural networks can automatically determine their tensor ranks in the training process.


\end{abstract}


\section{Introduction} 

As an efficient tool to overcome the curse of dimensionality,   
tensor decomposition methods date back to 1927 \cite{hitchcock1927expression} and have been employed in many application fields such as computer vision \cite{liu2012tensor}, signal processing \cite{cichocki2015tensor, sidiropoulos2017tensor}, graph matching \cite{cui2018quadratic},  bio-informatics \cite{durham2018predictd}, etc.  
Different from its matrix counterpart (i.e., singular value decomposition), tensor decompositions have different formats,  such as the CP decomposition \cite{hitchcock1927expression}, Tucker decomposition \cite{tucker1966some}, tensor-train decomposition \cite{oseledets2011tensor}, tensor network factorization~\cite{cichocki2016tensor, cichocki2017tensor}, t-SVD decomposition \cite{zhang2016exact}, and so forth. Some papers  have provided excellent surveys of tensor computation and its applications~\cite{kolda2009tensor, sidiropoulos2017tensor}. 

This paper will provide a high-level survey of tensor computation in the following two application fields: uncertainty-aware design automation and deep learning. These two seemingly irrelevant topics both require compact computational models to facilitate their subsequent statistical estimation, performance prediction and hardware implementation, despite their fundamentally different challenges:
\begin{itemize}[leftmargin=*]
    \item An EDA framework involves many modeling, simulation, and optimization modules. These modules often require some model-based simulation or hardware measurement data to decide the next step, but obtaining each piece of data is expensive. This challenge becomes more significant as process variations increase: one needs more data to capture an uncertain performance space. Therefore, it is desirable to extract high-quality ``compact" models to facilitate a decision making process with a ``small" available data set. 
    \item In deep learning, ``big" training data sets are often easy to obtain, and large-size neural networks can be trained on powerful platforms (e.g., in the cloud or on local high-performance servers). However, deploying them on resource-constrained hardware platforms (e.g., embedded systems and IoT devices) becomes a big challenge. As a result, there is a strong motivation to develop compact neural network models that can be deployed with low memory and computational cost. 
\end{itemize}
For both applications, tensor methods can be used to develop compact models with low computational and memory cost.

\section{Tensor Decomposition and Completion}

We first give a high-level tutorial about two important tensor problems: tensor decomposition and tensor completion. The first is often used to generate a compact low-rank representation when a (big) complete data set is given. The second is often employed when a (small) portion of the data is available. 

Fig.~\ref{fig:my_label} shows a tensor and several popular tensor decomposition techniques. A tensor $\ten A\in\mathbb{R}^{I_1\times I_2\times\ldots\times I_d}$ is a $d$-dimensional data array with $d\ge3$. It reduces to a matrix when $d=2$ and a vector when $d=1$. 
The high dimensionality of a tensor often brings in higher expressive power and higher compression capability. Three mainstream tensor decomposition methods are widely used for data analysis, scientific computing, and machine learning:
\begin{itemize}[leftmargin=*]
    \item The \textbf{CP decomposition}  \cite{hitchcock1927expression}  method decomposes a tensor into the summation of $R$ rank-one terms:
\begin{equation}\label{equ:CP}
\nonumber
    \ten A\approx\sum_{r=1}^R \mat x_r^{(1)} \circ \mat x_r^{(2)} \circ \ldots \circ \mat x_r^{(d)} 		\Leftrightarrow a_{i_1i_2 \cdots i_d}\approx\sum_{r=1}^R \prod_{k=1}^d\mat x_r^{(k)}(i_k).
\end{equation}
Here $\circ$ denotes the outer product and $\mat x_r^{(k)}\in\mathbb R^{I_i}$, and $a_{i_1i_2 \cdots i_d}$ denotes the scalar element in $\ten A$ indexed by $(i_1, i_2, \cdots, i_d)$. The total storage complexity is reduced to $\sum_{i=1}^dI_iR$. When the  approximation is replaced with  equality, the minimal integer of $R$ is called the {\sl tensor rank}. 

\item The \textbf{Tucker decomposition}  \cite{tucker1966some} compresses a tensor into a smaller core tensor $\ten G$ and $d$ orthogonal factor matrices $\{ \mat U_k \in \mathbb{R}^{I_k \times R_k} \}_{k=1}^d$: 
\begin{equation}
     a_{i_1\ldots i_d}=\sum_{r_1,\ldots,r_d=1}^{R_1,\ldots,R_d} g_{r_1\ldots r_d}\mat U_1(i_1,r_1)\ldots \mat U_d (i_d,r_d).
\end{equation}
The Tucker rank is bounded by $R_i\le I_i$ for all $i=1,\ldots,d$. The storage complexity is reduced to $\sum_{i=1}^dI_iR_i+\Pi_{i=1}^dR_i$. 

\begin{figure*}
    \centering
    \includegraphics[width=6in]{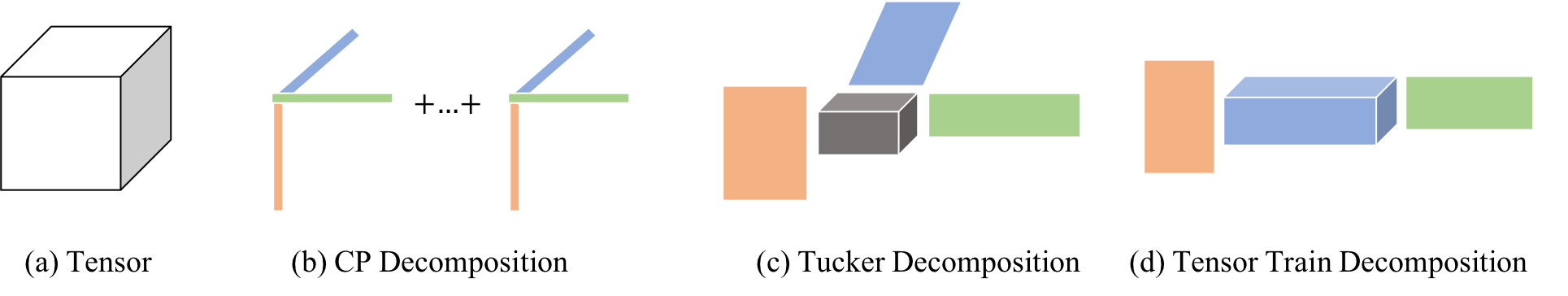}
    \caption{Several popular tensor compositions.}
    \label{fig:my_label}
\end{figure*}

\item The \textbf{tensor-train decomposition}  \cite{oseledets2011tensor} writes a tensor as a series of three-dimensional factor tensors, i.e., 
\begin{equation}
\label{eq: tt definition}
    a_{i_1 \ldots i_d}=\ten G_1 (:,i_1,:)\ten G_2 (:,i_2,:) \ldots \ten  G_d ({:,i_d,:}). 
\end{equation}
Here   $\ten G_k \in\mathbb{R}^{R_{k-1}\times I_k\times R_{k}}$, $R_0=R_d=1$, and $\ten G_k(:,i_k,:)$ is a MATLAB-like expression for the $i_k$-th lateral slice of $\ten G_k$. For a given tensor-train rank  $\mat{R}=\left( R_0,R_1,\dots, R_d\right)$, the storage complexity is reduced to  $\sum_{k=1}^dI_{k}R_{k-1}R_k$. 

\end{itemize}


\vspace{10pt}
{\bf Tensor Completion.} Given only partial elements of a tensor, the tensor completion or tensor recovery problem solves 
\begin{equation}\label{equ:tensorc}
    \min_{\ten X\in \ten M}\  \|P_{\Omega}(\ten A-\ten X)\|_F^2,
\end{equation}
where $\ten M$ denotes the set of low rank tensors in a proper format (e.g., CP,  Tucker or tensor-train format), the projection $P_{\Omega}(\ten A)$ keeps the element $a_{i_1 i_2 \cdots i_d}$ for all $(i_1,\ldots,i_d)\in\Omega$ and sets other elements to zero. The Frobenius norm is defined as $\|\ten A\|_F^2=\sum \limits_{i_1,\ldots,i_d}a_{i_1\ldots i_d}^2$. The cost function and regularization may be modified dependent on practical applications.
 
\section{Tensors for Uncertainty/Variability Analysis}

\subsection{Data-Expensive EDA Problems}
EDA problems are often model-driven and data-expensive: the design problems are well described by a detailed mathematical model (e.g., Maxwell equation for interconnect or RF device modeling, modified nodal analysis for circuit simulation), and one often needs to solve such an expensive mathematical model repeatedly or iteratively to get enough data (e.g., gradient information) to decide the next step (e.g., to optimize a circuit design parameter). The involved numerical computation makes the data acquisition expensive.  To accelerate the whole data-expensive EDA flow, one can choose to:
\begin{itemize}[leftmargin=*]
    \item Speed up the acquisition of each data sample. Representative matrix/vector-based techniques include fast PDE solvers~\cite{kamon1994fasthenry,nabors1991fastcap,phillips1997precorrected}, fast circuit simulators~\cite{kundert2013steady,li2012parallel} and model order reduction techniques~\cite{rewienski2003trajectory,daniel2004multiparameter,odabasioglu1997prima}. These solvers are often deterministic.
    \item Reduce the number of data acquisitions. This is especially important for nano-scale design that is highly influenced by process variations. In this case, one needs a huge amount of data samples to characterize the  uncertain circuit performance. Representative techniques include fast Monte Carlo~\cite{singhee2008practical} and recent stochastic spectral methods~\cite{zzhang:tcad2013, zzhang_cicc2014, manfredi:tcas2014}. 
\end{itemize}

 Tensor techniques can be employed to achieve both goals. Representative tensor techniques for the first goal include~\cite{liu2015model,liu2015staves,chen2019tensor}. In this paper, we focus on tensor techniques for uncertainty and variability analysis as summarized in Table~\ref{tab:my_label}, and we elaborate their key ideas below.

\subsection{Tensor Methods for Uncertainty Propagation}

Uncertainty quantification techniques aim to predict and control the probability density function (PDF) of the system output $y(\vecpar)$ under some random parameters $\vecpar\in\mathbb{R}^d$ describing process variations. The stochastic spectral methods based on generalized polynomial chaos method~\cite{gPC2002} have significantly outperformed Monte Carlo in many application domains. The key idea is to approximate $y(\vecpar)$ as a truncated linear combination of some specialized orthogonal basis functions $\{ \multiGPC_{\basisInd}(\vecpar) \}$ of $\vecpar$. The weight of each basis can be computed by various numerical techniques such as stochastic collocation \cite{xiu2005high}, stochastic Galerkin \cite{sfem} and stochastic testing \cite{zzhang:tcad2013}. When the parameters are non-Gaussian correlated, one can also employ the modified basis functions and stochastic collocation methods proposed in~\cite{Cui:EPEPS2018, cui2018stochastic,cui2018iccad,cui2019high}. 

\vspace{10pt}
{\bf Curse of Dimensionality.} Stochastic spectral methods suffer from an extremely high computational cost as the number of random parameters increase. For instance, in sampling-based techniques, the number of simulation samples may increase exponentially as $d$ increases. Some tensor solvers have been developed to address this fundamental challenge:
\begin{itemize}[leftmargin=*]
\item {\bf Tensor-Based Stochastic Collocation~\cite{zhang2016big,zhang2017big}.} The stochastic collocation method uses a projection method to compute the weight of each basis. Standard techniques discretize each random parameter into $m$ points, leading to $m^d$ simulation samples in total. Instead of simulating all samples, the technique in~\cite{zhang2017big} only simulates a small number of random samples and estimate the big unknown simulation data set by a tensor completion subject to two constraints: (1) the recovered tensor is low-rank; (2) the resulting generalized polynomial chaos expansion is sparse. This technique has been successfully applied to electronic IC, photonics and MEMS with up to $57$ random parameters. 


\item {\bf Tensor-Based Hierarchical Uncertainty Analysis~\cite{zhang2015enabling}.} Hierarchical techniques can be used to analyze the uncertainty of a complex system consisting of multiple interconnected components or subsystems. The key idea is to simulate each subsystem by a fast stochastic spectral method, then use their outputs as new random inputs for the system-level configuration. A major computational bottleneck is the high-dimensional integration required to recompute the basis function at the system level. In~\cite{zhang2015enabling} a tensor train decomposition is used to reduce the repetitive functional evaluation cost from an exponential cost to a linear one. This technique has enabled efficient uncertainty analysis of a MEMS/IC co-design with $184$ process variation parameters.

\item {\bf Tensor Method to Handle Non-Gaussian Correlated Uncertainties~\cite{cui2018iccad}.} A fundamental challenge of uncertainty propagation is how to handle non-Gaussian correlated process variations. Recently a set of basis functions and stochastic collocation methods were developed by~\cite{Cui:EPEPS2018, cui2018stochastic} to achieve high accuracy and efficiency. It is expensive to compute the basis functions in a high-dimensional case. When the non-Gaussian correlated random parameters are described by a Gaussian mixture density function, the basis functions were efficiently calculated by a functional tensor train method~\cite{cui2018iccad,cui2019high}. The key idea is as follows: the integration of a $d$-variable polynomial over each correlated Gaussian density can be written as the product of moments for each random variable. 
\end{itemize}

\begin{table*}[t]
    \centering
     \caption{Application of tensors in uncertainty propagation and variability Modeling}
    \begin{tabular}{c|c|c}
    \hline
   Reference & Problem & Key Idea  \\ \hline 
    \cite{zhang2017big,zhang2016big, konakli2016global, konakli2016reliability}&   high-dim stochastic collocation &   tensor completion to estimate unknown simulation data \\\hline
    \cite{zhang2015enabling}     & hierarchical uncertainty quantification & tensor-train decomposition for high-dim integration \\\hline
    \cite{cui2018iccad} & uncertainty analysis with non-Gaussian correlated uncertainty  &   functional tensor train to compute basis functions  \\\hline
    \cite{luan2019prediction} & spatial variation pattern prediction &   statistical tensor completion to predict variation pattern \\\hline
    
    \end{tabular}
   
    \label{tab:my_label}
\end{table*}



\subsection{Tensor Methods in Variability Prediction}

Statistical simulation of a circuit or device requires a given detailed statistical description (e.g., a probability density function) of the process variations. These statistical models are normally obtained by measuring the performance data of a huge number of testing chips. However, measuring the testing chips costs time and money, and may also cause mechanical damage. Existing techniques such as the virtual probe~\cite{zhang2011virtual} use the compressed sensing technique to predict the 2-D spatial variation pattern based on limited data.

In our recent paper~\cite{luan2019prediction}, we proposed to simultaneously predict the variation patterns of multiple dies. If each die has $N_1\times N_2$ devices to test, we can stack $N_3$ dies together to form a tensor. Then, we employed the Bayesian tensor completion technique in~\cite{zhao2015bayesian} to predict the spatial variation pattern with only a small number of testing samples. This technique can  automatically determine the tensor rank, achieving around $0.5\%$ relative errors with only $10\%$ testing samples with significant memory and computational cost reduction compared with the virtual probe technique~\cite{zhang2011virtual}.

\section{Compact Deep Learning Models}

Different from model-driven and data-expensive EDA problems, deep learning is suitable for data-driven and data-cheap applications such as computer vision and speech recognition. With the huge amount of available data (e.g., obtained from the social networks and many edge devices) and today's powerful computing platforms, deep learning has achieved success in a wide range of practical applications. However, deploying large neural networks requires huge computational and memory resource, limiting their applications on resource-constrained devices (i.e. smartphones, mobile robotics). Therefore, it is highly desirable  to build compact neural network models that can be deployed with low hardware cost. 


Many techniques can help generating more compact deep learning models. Most existing techniques are applied to individual weights, convolution filters or neurons, for instance: 
\begin{itemize}[leftmargin=*]
\item {\bf Pruning}~\cite{hanson1989comparing}: the key idea is to generate a sparse deep neural network by removing some redundant neurons that are not sensitive to the prediction performance. 

\item {\bf Quantization}~\cite{hinton1993keeping}: considering that model parameters are actually represented with binary bits in hardware, one may reduce the number of bits with little loss of accuracy. 

\item {\bf Knowledge Distillation}~\cite{hinton2015distilling}: the key idea is to shift the information from a deep and wide (teacher) neural network to a shallow one (i.e., a student network). 

\item {\bf Low-rank Compression}~\cite{denton2014exploiting}:  one can also compress the weight matrix or convolution filters by low-rank matrix or tensor decomposition.

\end{itemize} 

This section will survey the recent low-rank tensor techniques for generating compact deep neural network models. These techniques can be classified into two broad families:
\begin{itemize}[leftmargin=*]
    \item {\bf Tensorized Inference:} these techniques employ a ``train-then-compress" flow. Firstly a large deep neural network is trained (possibly with a GPU cluster), then tensor decomposition is applied to compress this pre-trained model   to enable its deployment on a hardware platform with limited resources (e.g., on a smart phone). 
    \item {\bf Tensorized Training:} these techniques skip the expensive training on a high-performance platform, and they aim to directly train a compact tensorized neural network from scratch and in an end-to-end manner.  
\end{itemize}

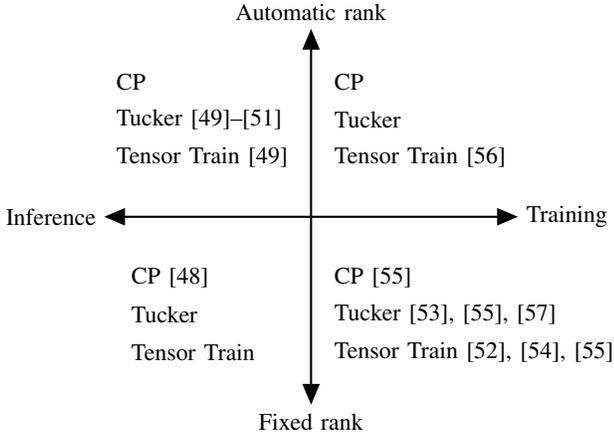
\begin{figure}

\small
\centering
 

\begin{tikzpicture}
\small 

\draw[thick, <->] (-2.75,0) -- (2.75,0);
\draw[thick, <->] (0,-2.5) -- (0,2.5);
\tkzText[above](0,2.5){Automatic rank};
\tkzText[below](0,-2.5){Fixed rank};
\tkzText[right](2.75,0){Training};
\tkzText[left](-2.75,0){Inference};

\tkzText[right](0.2,1.8){CP}
\tkzText[right](0.2,1.3){Tucker}
\tkzText[right](0.2,0.8){Tensor Train\cite{hawkins2019bayesian}}

\tkzText[right](-2.7,1.8){CP}
\tkzText[right](-2.7,1.3){Tucker \cite{kim2015compression, gusak2019one,garipov2016ultimate}}
\tkzText[right](-2.70,0.8){Tensor Train~\cite{garipov2016ultimate}}

\tkzText[right](0.2,-0.8){CP \cite{tjandra2018tensor}}
\tkzText[right](0.2,-1.3){Tucker \cite{kossaifi2017tensor,calvi2019tucker,tjandra2018tensor}}
\tkzText[right](0.2,-1.8){Tensor Train\cite{novikov2015tensorizing,tjandra2017compressing,tjandra2018tensor}}

\tkzText[right](-2.5,-0.8){CP \cite{lebedev2014speeding}}
\tkzText[right](-2.5,-1.3){Tucker}
\tkzText[right](-2.5,-1.8){Tensor Train}
    
\end{tikzpicture}
\caption{Summarize existing tensorized deep neural networks.}
\label{fig:class_DNN_compression}
\end{figure}

A common challenge of the above technique is to determine the tensor rank. Exactly determining a tensor rank in general is NP-hard \cite{haastad1990tensor}. Therefore, in practice one often leverages numerical optimization or statistical techniques to obtain a reasonable rank estimation. Dependent on the capability of automatic rank determination, existing tensorized neural network methods can be further classified into four groups as shown in Fig.~\ref{fig:class_DNN_compression}. We will elaborate their key ideas below.



\subsection{Tensorized Inference with a Fixed Rank} 
\label{subsec:tensorcompress}


Lebedev et al.~\cite{lebedev2014speeding} firstly applied CP tensor factorization to compress large-scale neural networks with fully connected layers. Because it is hard to automatically determine the exact CP tensor rank, this method keeps the tensor rank fixed in advance and employs an alternating least square method to compress folded weight matrices.



\subsection{Tensorized Inference with Automatic Rank Determination}

Compared with CP factorization, Tucker and tensor-train decompositions allow one to adjust the ranks based on accuracy requirement, and they have been employed in~\cite{garipov2016ultimate, kim2015compression} to compress both fully connected layers and convolution layers. A hardware prototype was even demonstrated for mobile applications in~\cite{kim2015compression}. Recently, an iterative compression technique was further developed to improve the compression ratio and model accuracy~\cite{gusak2019one}.

\subsection{Tensorized Training with a Fixed Rank}
In order to avoid the expensive pre-training in the uncompressed format, the work in \cite{novikov2015tensorizing} and \cite{calvi2019tucker} directly trained fully connected and convolution layers in low-rank tensor-train and Tucker format with the tensor ranks fixed in advance. This idea  has also been applied to recurrent neural networks \cite{tjandra2017compressing, tjandra2018tensor}.
 
\begin{figure*}[t]
    \centering
    \includegraphics[width=5in]{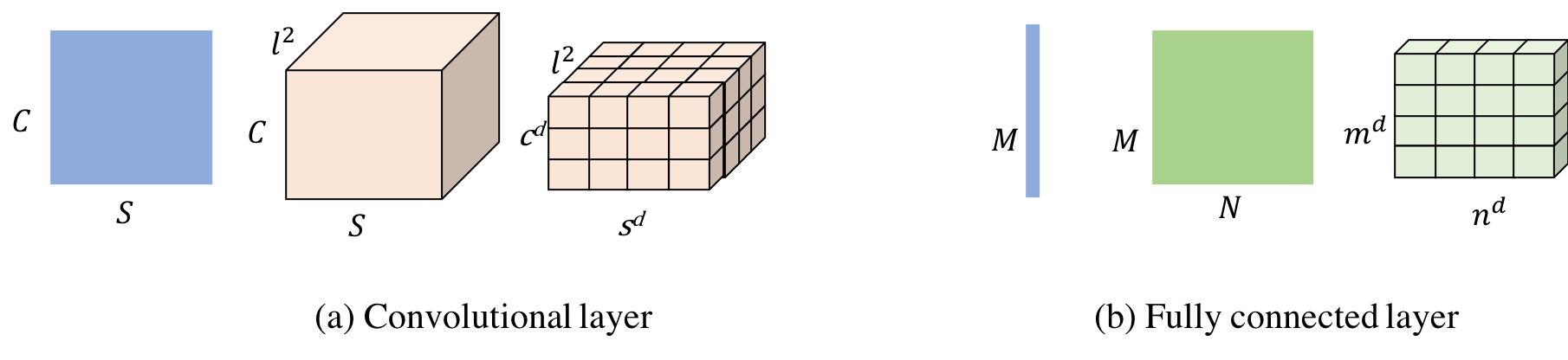}
    \caption{(a) The convolutional kernel can be regarded as a 3D tensor in $C\times S\times l^2$ or a 4D tensor in $C\times S\times l\times l$. We can also reformulate it into a ($2d+2$)-dimensional tensor. (b) The weight matrix in the fully connected layer is a 2D matrix, which can be reformulated into a $2d$-dimensional tensor.
    }
    \label{fig:tensor_DNN}
\end{figure*}

In the following, we would like to summarize the storage and computational complexity (e.g. flops) by applying CP decomposition, Tucker decomposition, and tensor-train decomposition for both the convolution layers and the fully connected layers: 
\begin{itemize}[leftmargin=*]
    \item {\bf Tensorized Convolution Layers.} Consider a convolutional  weight $\ten{K}\in \mathbb{R}^{l\times l \times C \times S }$, where $l$ is the filter size, and $C, S$ are the number of input and output channels respectively. This kernel requires $l^2CS$ parameters in total. The  tensor-train format  reduces the number of parameters from $l^2CS$ to $O(R^2(2l+C+S))$. The rank-$R$ CP-decomposition can represent the weight with $O(R(2l+C+S))$ storage.   
The Tucker-2 decomposition  \cite{kossaifi2017tensor} only compresses the input channel and output channel dimension, and requires $O(R(Rl^2+C+S))$ storage. 
We also note that in practice it is common to reshape the tensor $\ten{K}$ into a higher-order tensor [as shown in Fig.~\ref{fig:tensor_DNN}~(a)]. Specifically, one can factorize 
\begin{equation*}
  C=\prod_{k=1}^d c_k \text{ and } S=\prod_{k=1}^d s_k
\end{equation*}
and reformulate $\ten K$ as a $(2d+2)$-dimensional tensor $\ten K\in\mathbb{R}^{l\times l\times c_1\ldots\times c_d\times s_1\ldots\times s_d}$.  
 Such folding can often result in a higher compression ratio.

\item {\bf Tensorized Fully Connected Layers.} The weight in the fully connected layer is a 2D-matrix, and direct application of the tensor method reduces to a matrix singular value decomposition.  Given a weight matrix $\mat{W}\in \mathbb{R} ^{M\times N}$ one can factorize $M=\prod_{k=1}^d M_k$,  $N=\prod_{k=1}^d N_k$ 
and represent $\mat{W}$ with a high-dimensional tensor $\ten{W}\in\mathbb{R}^{m_1\times n_1\times\ldots\times m_d\times n_d}$ [as shown in Fig.~\ref{fig:tensor_DNN}~(b)]. 
Then by the tensor factorization, the number of parameters can reduce  from $M\times N$ to $O(r(md+nd))$,   and $O(dr^2mn)$ for a CP format  and tensor-train format, respectively. 
\end{itemize}


We summarize both the storage and computational complexity of different tensor compression methods in Table~\ref{tab:DNN_complexity}. For the convolutional layer, we only counts the computational costs of a $l\times l\times C$ block. 

\begin{table*}[]
\caption{Speedup of storage and computational  complexity of tensor methods for the  convolutional (conv) layer and  the fully connected (FC) layer. `-' denotes inapplicable: Tucker decomposition is not applicable to high-order tensors due to the curse of dimensionality; FC weight is a 2D matrix.}
    \centering
     
    \begin{tabular}{c|cc|c|cc|cc}
    \hline
    
      &    \multicolumn{2}{c|}{ Original} & & \multicolumn{2}{c|}{ Tensor decomposition} & \multicolumn{2}{c}{High-order tensor decomposition} \\ 
      &  Storage & FLOPS & & Storage & FLOPS & Storage & FLOPS \\\hline
      &   && CP & $O(R(2l+C+S))$ & $O(R(Cl^2+S))$ &$O(r(2l+cd+sd))$ & $O(r(Cl^2+S))$ \\
     Conv & $O(CSl^2)$ & $O(CSl^2)$&  Tucker &$O(R(Rl^2+C+S))$  & $O(R(Cl^2+Rl^2+S))$ &-& -\\
     & &&TT & $O(R^2(2l+C+S))$  & $O(R^2(2Cl^2+C+S))$& $O(r^2(2l+dc+ds))$ & $O(r^2(2Cl^2+dC+dS))$\\ 
     \hline 
     
     &   && CP & $O(R(C+S))$ & $O(R(C+S))$&$O(dmr+dnr)$ & $O(Mr+Nr)$  \\
   FC & $O(MN)$ & $O(MN)$ &  Tucker & - &-& -& - \\
     & &&TT &- &-& $O(r^2dmn)$ & $O(r^2md \max(M,N))$ \\ 
     \hline 
   
    \end{tabular}
      
    \label{tab:DNN_complexity}
\end{table*}


\subsection{Tensorized Training with Automatic Rank Determination} 
It is a challenging task to automatically determine the tensor rank in a training process due to the following reasons:
\begin{itemize}[leftmargin=*]
    \item Different from the matrix case, there is not a proper surrogate model for the rank of a high-order tensor. As a result, it is non-trivial to regularize the loss function of a neural network with a low-rank penalty term.
    \item Existing tensor decomposition methods work on a given full tensor. However, the tensors in   end-to-end training are embedded within a deep neural network in a highly nonlinear manner. This also makes existing tensor completion or recovery frameworks fail to work.
\end{itemize}

In order to address this fundamental challenge and to enable efficient end-to-end training, we leveraged the variational inference and proposed the first Bayesian tensorized neural network~\cite{hawkins2019bayesian} to automatically determine the tensor rank as part of the training process. 
We use the notation $\mathcal{D} = \{(\mat{x}_i,\mat{y}_i)\}_{i=1}^N$ to represent the training data. 
Our goal is to learn the low-rank tensor parameters by estimating the following posterior density: 
\begin{equation}
   p(\boldsymbol{\theta} | \mathcal{D}) \propto p( \mathcal{D} | \boldsymbol{\theta} ) p(\boldsymbol{\theta})
\end{equation} 
where $\boldsymbol{\theta} $ include the unknown tensor-train factors and some hyper-parameters, and $p(\boldsymbol{\theta})$ is a prior density to enforce low-rank property. The following key techniques enabled the efficient training and automatic rank determination:
\begin{itemize}[leftmargin=*]
    \item The prior density $p(\boldsymbol{\theta})$ is designed by considering the coupling of adjacent tensor-train cores, such that their ranks can be controlled simultaneously.  
    \item We employed a Stein variational gradient descent~\cite{liu2016stein} method to approximate the posterior density $ p(\boldsymbol{\theta} | \mathcal{D})$. This method combines the flexibility of Markov-Chain Monte Carlo and the efficiency of optimization techniques, which is beyond the capability of the mean-field inference framework in the Bayesian tensor completion framework~\cite{zhao2015bayesian}.
\end{itemize}

This method has trained a two-layer fully connected neural network, a 6-layer CNN and a 110-layer residual neural network, leading to $7.4\times$ to $137\times$ compression ratios.

\section{Conclusions}
In this paper, we have revisited several compact models generated by the tensor decomposition/completion approach. 
For the data-expansive problems arising from EDA, we have summarized several tensor methods in uncertainty quantification and spatial prediction. Tensor techniques have successfully solved many high-dimensional uncertainty quantification problems with both independent and non-Gaussian random parameters. They have also significantly reduced the chip testing cost in spatial variation pattern prediction.

In the context of deep learning, tensor decomposition proves to be an efficient technique to obtain compact learning models. They have achieved significant compression in both inference and training. Our recent Bayesian tensorized neural network allows automatic tensor rank determination in the end-to-end training process.

\bibliographystyle{IEEEtran}
\bibliography{ref}
\end{document}